\documentclass{mcom-l}
\usepackage[utf8]{inputenc}
\usepackage[T1]{fontenc}
\usepackage{ae,aecompl}
\usepackage{amsmath, amssymb}
\usepackage{url}

\newcommand{\OD}{\mathcal O_D}
\newcommand{\Sl}{\mathrm {Sl}_2 (\Z)}

\newcommand{\N}{\mathbb N}
\newcommand{\Z}{\mathbb Z}
\newcommand{\Q}{\mathbb Q}
\newcommand{\C}{\mathbb C}
\newcommand{\F}{\mathbb F}
\newcommand{\Hb}{\mathbb H}
\newcommand{\Fc}{\mathcal F}
\newcommand{\f}{\mathfrak f}
\newcommand{\pf}{\mathfrak p}
\newcommand{\w}{\mathfrak w}
\newcommand{\wN}{\mathfrak w_N}
\newcommand{\wpq}{\mathfrak w_{p_1, p_2}}

\newtheorem{theorem}{Theorem}[section]
\newtheorem{corollary}[theorem]{Corollary}
\newtheorem{algorithm}[theorem]{Algorithm}
\newcommand{\For}{{\textbf {for} }}
\newcommand{\To}{{\textbf {to} }}
\newcommand{\Downto}{{\textbf {downto} }}
\newcommand{\Return}{{\textbf {return} }}

\copyrightinfo{2008}{Andreas Enge}

\numberwithin{equation}{section}

\begin{document}

\title [Complexity of floating point class polynomial computation]
{The complexity of class polynomial computation via
floating point approximations}

\author {Andreas Enge}
\address {INRIA Saclay--\^Ile-de-France
\& Laboratoire d'Informatique (CNRS/UMR 7161),
\'Ecole polytechnique, 91128 Palaiseau Cedex, France}
\email {enge@lix.polytechnique.fr}

\date {April 23, 2007; revised April 30, 2008}


\begin {abstract}
\footnotesize
We analyse the complexity of computing class polynomials, that are an important
ingredient for CM constructions of elliptic curves, via complex floating point
approximations of their roots. The heart of the algorithm is the evaluation
of modular functions in several arguments. The fastest one of the presented
approaches uses a technique devised by Dupont to evaluate modular functions
by Newton iterations on an expression involving the arithmetic-geometric mean.
Under the heuristic assumption, justified by experiments, that the correctness of the result is not perturbed by rounding errors, the algorithm runs in time
\[
O \left( \sqrt {|D|} \log^3 |D| \, M \left( \sqrt {|D|} \log^2 |D| \right) \right)
\subseteq O \left(|D| \log^{6 + \varepsilon} |D| \right)
\subseteq O \left( h^{2 + \varepsilon} \right)
\]
for any $\varepsilon > 0$, where $D$ is the CM discriminant, $h$ is the
degree of the class polynomial and $M (n)$ is the time needed to multiply two $n$-bit numbers.
Up to logarithmic factors, this running time matches the size of the constructed
polynomials. The estimate also relies on a new result concerning the complexity
of enumerating the class group of an imaginary quadratic order and on a
rigorously proven upper bound for the height of class polynomials.
\end {abstract}

\maketitle

\section {Motivation and results}

The theory of complex multiplication yields an efficient approach to the
construction of elliptic curves over a finite field having a given
endomorphism ring of discriminant $D$, as long as $|D|$ is
not too large. Exploiting the link between the endomorphism ring of an
elliptic curve and its number of points, it is possible to efficiently obtain
curves with specific properties. Applications include primality
proving \cite {am93}, the construction of classical, discrete logarithm
based elliptic curve cryptosystems and of identity based cryptosystems
\cite {mnt01,dem05,bls03con,bw05}.

The classical approach to effective complex multiplication is to compute
minimal polynomials of special, algebraic values of modular
functions on the upper complex halfplane (more details are provided in
Section~\ref {sec:cm}). The complexity of this algorithm remains a shady
issue; one of the reasons a serious analysis has not been undertaken so far
are the purported numerical difficulties with the computations with complex
floating point numbers, an issue that actually does not arise in practice
(cf. the short discussion at the end of Section~\ref {sec:implementation}).

In \cite {ch02}
the authors present an algebraic generalisation of the complex
multiplication algorithms to $p$-adic fields, that yields the same minimal
polynomials. Besides ruling out any possible numerical instabilities, they
obtain a complexity of $O (|D|^{1 + \varepsilon})$, which is better than a
straightforward implementation of the complex approach, and asymptotically
optimal (up to logarithmic factors) due to the size of the constructed objects.

The present article thus has two goals. First, it provides an accurate account
of the complexity of different algorithms for class polynomial computation via
floating point approximations. Second, it shows that asymptotically optimal
algorithms exist also in the complex setting, with a complexity that is linear
(up to logarithmic factors) in the output size. These new algorithms are
presented in Sections~\ref {ssec:multi} and \ref {ssec:agm}, and the faster one
achieves the following:

\begin {theorem}
\label {th:time}
Let $f$ be a fixed modular function that is a class invariant for a family of
discriminants $D$ of class numbers $h = h (D)$. Then the algorithm of
Section~\ref {ssec:agm}, which computes a floating point approximation to the
class polynomial for $f$, runs in time
\[
O \left( h (\log^2 h + \log n) \, M (n) \right)
\]
when executed with complex floating point numbers of $n = n (D)$ bits precision,
where $M (n)$ is the time needed to multiply two such numbers as detailed in Section~\ref {ssec:multiprec}.
\end {theorem}

The floating point precision $n$ required to carry out the computations is
clearly bounded from below by the height of the class polynomial. The following
theorem provides a rigorously proven upper bound on the heights, that is close
to experimental findings. A bound of the same shape is given in
\cite [\S 3.3]{alv04} as a heuristic, and the gist of the proof appears already
in \cite [\S 5.10]{ll90} and \cite [Section~2]{sch91exp}.

\begin {theorem}
\label {th:height}
The logarithmic height of the class polynomial for $j$ for the discriminant
$D$ of class number $h$ is bounded above by
\begin {eqnarray*}
&&
c_5 h + c_1 N \left( \log^2 N + 4 \gamma \log N
+ c_6 + \frac {\log N + \gamma + 1}{N} \right) \\
& \leq &
c_1 N \log^2 N + c_2 N \log N + c_3 N + c_1 \log N + c_4
\end {eqnarray*}
where $N = \sqrt {\frac {|D|}{3}}$,
$\gamma = 0.577\ldots$ is Euler's constant,
$c_1 = \sqrt 3 \pi = 5.441\ldots$,
$c_2 = 18.587\ldots$,
$c_3 = 17.442\ldots$,
$c_4 = 11.594\ldots$,
$c_5 = 3.011\ldots$ and
$c_6 = 2.566\ldots$.

\noindent
The asymptotic upper bound of
\[
O \left( \sqrt {|D|} \log^2 |D| \right)
\]
holds for any other class invariant as well.
\end {theorem}

In practice, one observes that rounding errors do not disturb the result. It
suffices to take $n$ as an approximation of the height plus a few guard digits
to be able to round the floating point approximation of the class polynomial to
the correct polynomial with integral coefficients. A rigorous error analysis,
however, appears to be out of reach. So Theorems~\ref {th:time} and \ref
{th:height} can be brought together only in the form of a heuristic, assuming
that the computations with floating point numbers of $n$ bits yield an
approximation of the class polynomial that is correct on essentially $n$ bits.

\begin {corollary}[heuristic]
\label {cor:time}
Taking $n \in O \left( \sqrt {|D|} \log^2 |D| \right)$ in Theorem~\ref {th:time}
and using the bound on the class number $h \in O \left( \sqrt {|D|} \log |D|
\right)$ proved at the end of
Section~\ref {sec:height}, the algorithm of Section~\ref {ssec:agm} computes the
class polynomial for $f$ in time
\[
O \left( \sqrt {|D|} \log^3 |D| \, M \left( \sqrt {|D|} \log^2 |D| \right) \right)
\subseteq O \left( |D| \log^{6 + \varepsilon} |D| \right)
\subseteq O \left( h^{2 + \varepsilon} \right)
\]
for any $\varepsilon > 0$.
\end {corollary}

Notice that up to logarithmic factors, this complexity corresponds to the output
size of the algorithm, namely the size of the class polynomials.
Notice also that the correctness of the output can be verified by a
probabilistic, Monte-Carlo type algorithm; namely one may check that the
reductions of the class polynomial modulo sufficiently many suitable primes $p$
yield elliptic curves over $\F_p$ with complex multiplication by $\OD$. Indeed,
the main application of class polynomials is to compute elliptic curves over
finite fields with given complex multiplication, and in this situation
it can be verified independently that the curves are correct.

As an ingredient for the proof of Theorem~\ref {th:time} we obtain in
Section~\ref {sec:class group} the following result for computing class groups
of imaginary quadratic orders:

\begin {theorem}
\label {th:classgroup}
The class group of the imaginary-quadratic order of discriminant $D$
and class number $h$ can be enumerated
\begin {itemize}
\item
unconditionally by a probabilistic algorithm in time
\[
O \left( \sqrt {|D|} \log |D| \log\log |D| M (\log |D|) \right);
\]
\item
under GRH by a deterministic algorithm in time
\[
O \left( h \log\log |D| M (\log |D|) \right).
\]
\end {itemize}
\end {theorem}
Again the algorithm is essentially optimal in its output size.

\section {Complex multiplication and class polynomials}
\label {sec:cm}

\subsection {The basic approach}

For proofs of the following facts on complex multiplication of elliptic
curves see, for instance, \cite {cox89}.

Let us first consider the situation over the complex numbers. Let $D < 0$
be an imaginary quadratic discriminant, and $\OD = \left[ 1, \frac { D +
\sqrt D}{2} \right]_\Z$ the (not necessarily maximal) order of discriminant
$D$ in $K = \Q (\sqrt D)$. The ideal class number of $\OD$ is denoted by $h
= h_D$. By Siegel's theorem \cite {sie36},
$\frac {\log h}{\log |D|} \to \frac {1}{2}$ $(|D| \to \infty)$, so that
$|D| \in O (h^ {2 + \varepsilon})$ and $h \in O (|D|^{1/2 + \varepsilon})$
for any $\varepsilon > 0$. There are $h$ isomorphism classes of elliptic
curves over $\C$ having {\em complex multiplication} by $\OD$, that is,
curves with $\OD$ as endomorphism ring. Namely, let
$j : \Hb = \{ z \in \C : \Im (z) > 0 \} \to \C$
denote the absolute modular invariant, and let
$\tau_i = \frac {-B_i + \sqrt D}{2 A_i}$
run through the roots in $\Hb$ of the reduced quadratic forms
$[A_i, B_i, C_i] = A_i X^2 + B_i X + C_i$ of discriminant
$D = B_i^2 - 4 A_i C_i$,
representing the ideal classes of $\OD$; then the
j-invariants of the elliptic curves are given by the $j (\tau_i)$.
Moreover, these $j (\tau_i)$ are algebraic integers; in fact, they generate
the so-called {\em ring class field} $K_D$ for $\OD$, the Galois
extension of $K$ whose Galois group is isomorphic to the class group of $\OD$
(the isomorphism being given by the Artin map).
The minimal polynomial of the $j (\tau_i)$ over $K$,
$H_D (X) = \prod_{i=1}^h (X - j (\tau_i))$,
has in fact coefficients in $\Z$ and is called a {\em class polynomial}. In
the special case that $D$ is a fundamental discriminant, the ring class field
$K_D$ is also called the {\em Hilbert class field} of $K$.

Let now $\F_q = \F_{p^m}$ be a finite field of characteristic $p$. Suppose
that $p$ splits in $K = \Q (\sqrt D)$ and that $p \nmid D$; then
$p$ is unramified in $K_D$. If furthermore $q$ may be written as
$4 q = U^2 + D V^2$ with $U$, $V \in \Z$, then the inertia degree of
the prime ideals above $p$ in $K_D$ divides $m$. The reductions of the complex
elliptic curves with complex multiplication by $\OD$ modulo any of these
prime ideals thus live in $\F_q$. By Deuring's reduction and lifting theorems
\cite [Einleitung, par.~5]{deu41typ}, these $h$ curves are precisely the
elliptic curves over $\F_q$
with complex multiplication by $\OD$. They may be obtained as follows:
Compute the class polynomial $H_D \in \Z [X]$ and reduce it modulo $p$. It
splits completely over $\F_q$, and each of its roots is the j-invariant of
an elliptic curve over $\F_q$ with the desired endomorphism ring.

\subsection {Class invariants}

Unfortunately, $H_D$ has very large coefficients (see the discussion in
Section~\ref {sec:height}), so that its computation requires a high precision
to be accurate. In practice, one may often gain a constant factor for the
required number of digits by using instead of $j$ modular functions $f$
that are invariant under
$\Gamma^0 (N) = \left\{
\begin {pmatrix} a & b \\ c & d \end {pmatrix}
\in \Sl : N | b \right\}$
for some positive integer $N$; that is,
$f \left( \frac {a z + b}{c z + d} \right) = f (z)$
for any such matrix
$\begin {pmatrix} a & b \\ c & d \end {pmatrix}$.
Under suitable conditions on the discriminant $D$ and suitable
normalisations of the $\tau_i$,
derived from Shimura's reciprocity law, the {\em singular values} $f
(\tau_i)$ are still elements of the class field $K_D$ \cite
{sch76,gee99,gs98,sch02}; we then call $f$ a {\em class invariant} and the
minimal polynomial
$H_D [f] (X) = \prod_{i=1}^h (X - f (\tau_i))$ again a {\em class
polynomial}.
The first such class invariants are given by Weber's functions $\f$, $\f_1$,
$\f_2$, $\gamma_2$ and~$\gamma_3$ \cite {web08}.

Two parameterised families of class invariants are exhibited in \cite
{es04}, where $N$ is the product of two primes and $H_D [f] \in \Z [X]$,
and \cite {em07}, where $N$ is prime and $H_D [f]$ has coefficients in the
maximal order of $\Q (\sqrt D)$.
Again, these class polynomials split completely after reduction in $\F_q$.
The corresponding elliptic curves may be recovered from some polynomial
relationship
between $f$ and $j$: If the modular curve $X_0 (N)$ has genus $0$, then $j$
is sometimes given by a rational formula in $f$; otherwise it has been
suggested in \cite {es04} to look for a root of the modular polynomial
$\Phi (f, j)$ after reducing modulo $p$ and specialising in the value found
for $f$ in $\F_q$. (Both cases require that all coefficients be rational
to make sense modulo $p$, which holds for all exhibited class invariants.)

\section {Complexity of arithmetics}

This section discusses the well-known complexity of the basic multiprecision and
polynomial arithmetic underlying the computations.

\subsection {Multiprecision floating point arithmetic}
\label {ssec:multiprec}

Let $M (n)$ be the bit complexity of multiplying two $n$-bit integers. Using Schön\-hage--Stras\-sen multiplication \cite{ss71}, one has $M (n) \in O (n \log n \log \log n)$. With F\"urer's algorithm \cite{fue07}, one has $M (n) \in n \log n \, 2^{O (\log^* n)}$, where $\log^* n$ is the number of times the logarithm function has to be applied to $n$ before the result drops below $1$. So with either algorithm, $M (n) \in O \left( n \log^{1 + \varepsilon} n \right)$ for any $\varepsilon > 0$.

The four basic arithmetic operations and
the square roots of real floating point numbers of precision $n$ have a
complexity of $O (M (n))$, the inversion and square roots being realised by
Newton iterations as explained in Lemmata~2.2 and 2.3 of \cite {bre76}. The
same article shows that $\exp$, $\sin$ and the constant $\pi$ can be computed in
$O (\log n \, M (n))$.

Thus, the four basic operations on complex floating point numbers of
precision~$n$ can be executed in time $O (M (n))$. Letting
$c = \sqrt {\frac {a + \sqrt {a^2 + b^2}}{2}}$ for $a > 0$, one obtains
$\sqrt {a + b i} = c + \frac {b}{2c}$ in time $O (M (n))$. Finally, the complex
exponential may be reduced to the real exponential and the real sine and cosine
functions and is thus computed with complexity $O (\log n \, M (n))$.

\subsection {Polynomial arithmetic}
\label {ssec:poly}

Concerning operations with polynomials, we assume that a floating point
precision of $n$ bits has been fixed and is the same for all input and output
polynomials.
Multiplying two polynomials of degree $d$ over $\C$ by the FFT takes $O (d \log
d)$ multiplications in $\C$, whence it has a complexity of $O (d \log d \, M (n))$
once the necessary roots of unity have been computed.

The following algorithm (\cite[Algorithm 10.3]{gg99}) obtains a monic
polynomial of degree $h$ from its roots by organising the computations in a
binary tree.

\begin {algorithm}[\texttt {Poly\_from\_roots}]
\label {alg:poly}
\quad \\
\noindent
\textbf {Input:} $h \in \N$ and $x_1, \ldots, x_h \in \C$

\noindent
\textbf {Output:} binary tree
$T_{k,i} = \left( X - x_{(i-1) 2^k + 1} \right) \cdots
\left( X - x_{\min (i \, 2^k, h)} \right)$ \\
\hspace* {5mm}
for $k = 0, \ldots, t = \lceil \log_2 h \rceil$,
$i = 1, \ldots, 2^{t-k}$; \\
\hspace* {5mm}
in particular, the root
$T_{t,1} = (X - x_1) \cdots (X - x_h)$

\begin {enumerate}
\item
\For $i=1$ \To $h$ \\
\hspace* {5mm} $T_{0,i} \leftarrow X - x_i$
\item
\For $i=h+1$ \To $2^t$ \\
\hspace* {5mm} $T_{0,i} \leftarrow 1$
\item
\For $k=1$ \To $t$ \\
\hspace* {5mm} \For $i=1$ \To $2^{t-k}$ \\
\hspace* {10mm} $T_{k,i} \leftarrow T_{k-1,2i-1} \cdot T_{k-1,2i}$
\item
\Return $T$
\end {enumerate}
\end {algorithm}

If the needed roots of unity are precomputed, the algorithm has a complexity of
$O (h \log^2 h \, M (n))$. Notice that the roots of unity of order
$2^{\lceil \log_2 h \rceil}$ suffice to carry out all the FFTs in
Algorithm~\ref {alg:poly}. By Section~\ref {ssec:multiprec}, computing a
primitive root of unity takes
$O (\log n \, M (n))$, all others can be obtained by successive multiplications
in $O (h \, M (n))$. Thus, the total complexity of Algorithm~\ref {alg:poly}
becomes
\[
O ((h \log^2 h + \log n) M (n)).
\]

One of the asymptotically optimal algorithms for class polynomials relies on a
related technique of symbolic computation (see \cite [\S 10.1]{gg99}). Let
$f (X)$ be a polynomial of degree $d$
over $\C$ that is to be evaluated in $h$ different arguments
$x_1, \ldots, x_h$. The key observation for doing so fast is that $f (x_i)$ is
nothing but $f (X) \bmod X - x_i$. So the first step of the algorithm is
the construction of the binary tree $T_{k,i}$ as in Algorithm~\ref
{alg:poly}, containing products of more and more of the $X - x_i$. Then a
matching tree $R_{k,i} = f \bmod T_{k,i}$ is computed in the converse
order, from its root to its leaves, such that the leaves contain the
desired values of $f$.

\begin {algorithm}[\texttt {Multi\_eval}]
\label {alg:multi}
\quad \\
\noindent
\textbf {Input:} $h \in \N$, $f \in \C [X]$ and $x_1, \ldots, x_h \in \C$

\noindent
\textbf {Output:} $f (x_1), \ldots, f (x_h)$
\begin {enumerate}
\item
$T \leftarrow \text {{\tt Poly\_from\_roots} } (h, x_1, \ldots, x_h)$
\item
$t \leftarrow \lceil \log_2 h \rceil$, $R_{t,1} \leftarrow f \bmod T_{t,1}$
\item
\For $k=t-1$ \Downto $0$ \\
\hspace* {5mm} \For $i=1$ \To $2^{t-k}$ \\
\hspace* {10mm} $R_{k,i} \leftarrow R_{k+1, \left\lfloor \frac {i+1}{2}
\right\rfloor} \bmod T_{k,i}$
\item
\Return $R_{0,1}, \ldots, R_{0,h}$
\end {enumerate}
\end {algorithm}

A division with remainder of a polynomial of degree $d$ by a polynomial of
degree $d'$ has the same complexity as multiplying polynomials of degree bounded
by $d$ (\cite [\S 9.1]{gg99}). The
algorithm first computes an inverse modulo some power of $X$ by Newton
iterations and obtains the quotient with one multiplication; the remainder then
requires another multiplication. So after the first reduction of $f$ modulo
the product of all the $X - x_i$, the complexity of Algorithm~\ref
{alg:multi} on $n$-bit numbers is the same as that of
Algorithm~\ref {alg:poly}. Including the first reduction, it is given by
\[
O ((d \log d + h \log^2 h + \log n) M (n)).
\]

\section {The height of class polynomials}
\label {sec:height}

The running time of the algorithms working with floating point approximations
depends crucially on the precision $n$. To be able to round the approximated
class polynomial to a polynomial over the integers ($\Z$ or the maximal order of
an imaginary-quadratic number field), $n$ has to be at least the bit size of
the largest coefficient, or otherwise said, the logarithmic height of the
polynomial. In this section, we prove Theorem~\ref {th:height} by developing an
explicit upper bound without hidden constants for the height of the class
polynomial for $j$. The result is of independent interest; for instance, it
allows to bound the precision also for the $p$-adic algorithms, which yields a
proof of correctness for their output.

As shown in \cite {em02}, the height of class polynomials for a different $f$
changes asymptotically by a constant factor, depending on the degrees in $f$ and
$j$ of the modular polynomial connecting the two functions. So this section
proves an asymptotic bound of $O \left( \sqrt {|D|} \log^2 |D| \right)$ for the height of any
class polynomial. To obtain a more explicit bound for invariants other than $j$,
the techniques of this section may be used with the necessary adaptations.

Let us first give the intuitive basis for the bound.
Since the values of $j$ are usually large, more often than not the largest
coefficient of the class polynomial is its constant one. Then approximating
$j (\tau) = q^{-1} + 744 + \sum_{\nu = 1}^\infty c_\nu q^\nu$ by the first term,
$q^{-1}$, one obtains a heuristic estimate for the height as
\[
\pi \sqrt {|D|} \sum_{[A, B, C]} \frac {1}{A},
\]
where the sum is taken over all reduced primitive quadratic forms of
discriminant $D$. This heuristic estimate
is very well confirmed by experimental findings, see \cite {em02}.

Turning this idea into an explicit bound is mainly a matter of computations that
are sketched in the following.

The reducedness of the quadratic forms is equivalent to the $\tau$-values lying in the standard fundamental domain $\Fc$ for the action of $\Sl$ on $\Hb$,
\begin {equation}
\label {eq:fundamental domain}
\Fc = \left\{ z \in \Hb : |z| > 1, - \frac {1}{2} \leq \Re (z)
< \frac {1}{2} \right\}
\cup
\left\{ z \in \Hb : |z| = 1, - \frac {1}{2} \leq \Re (z)
\leq 0 \right\}.
\end {equation}
Then $|q| \leq \left| e^{2 \pi i \frac {1 + \sqrt {-3}}{2}} \right|
= e^{- \pi \sqrt 3}$.
The upper bound $c_\nu \leq \frac {e^{4 \pi \sqrt \nu}}{\sqrt 2 \nu^{3/4}}$
of \cite {bp05} then yields
\[
| j (\tau) - q^{-1} |
\leq 744 + \sum_{\nu = 1}^\infty \frac {e^{4 \pi \sqrt \nu}}{\sqrt 2 \nu^{3/4}}
e^{- \pi \sqrt 3 \, \nu}
= k_1 = 2114.566\ldots
\]
for $\tau$ in the fundamental domain $\Fc$, so that
\[
|j (\tau)|
\leq |q^{-1}| + k_1
\leq k_2 |q^{-1}|
\]
with
$k_2 = 1 + k_1 e^{-\pi \sqrt 3} = 10.163\ldots$

Let the $A_i$ be numbered in increasing order. Then the logarithm of the
absolute value of the coefficient in front of $X^j$ is bounded above by
\[
\log \left( \binom {h}{j} \prod_{i=1}^{h-j} (k_2 |q_i^{-1}|) \right)
\leq \log \left( 2^h k_2^h \prod_{i=1}^h |q_i^{-1}| \right)
\leq \log (2 k_2) h + \pi \sqrt {|D|} \sum_{i=1}^h \frac {1}{A_i},
\]
independently of $j$.

The next step consists of estimating $\sum_{i=1}^h \frac {1}{A_i}$.
It is proved in \cite [Lemma~2.2]{sch91exp} that
$\sum_{i=1}^h \frac {1}{A_i} \in O (\log^2 |D|)$.
We shall derive a bound that makes the involved constants explicit.
Consider the number of possible $B$ for a given $A$.
This number is certainly bounded above by the number of $B \in ]-A, A]$
that satisfy $B = D \bmod 2$ and $\frac {B^2 - D}{4} = 0 \bmod A$. Assuming the
worst case that the quadratic equation has a root modulo each prime divisor
of $A$ and considering the cases of odd and even $A$ separately, one obtains
an upper bound on the number of $B$ of $2 \cdot 2^{\omega (A)}
\leq 2 \tau (A)$, where $\omega (A)$ denotes the number of prime factors of
$A$ and $\tau (A)$ its number of divisors. Hence,
\[
\sum_{i=1}^h \frac {1}{A_i} \leq 2 \sum_{A=1}^{\sqrt {\frac {|D|}{3}}}
\frac {\tau (A)}{A},
\]
and this sum may be bounded by standard techniques from analytic number
theory. We use the estimates
\begin {equation}
\label {eq:harmonic}
\log N + \gamma \leq \sum_{n=1}^N \frac {1}{n} \leq
\log N + \gamma + \frac {1}{2N}
\end {equation}
with Euler's constant $\gamma = 0.577\ldots$ and
\begin {equation}
\label {eq:sumlog}
\sum_{n=1}^N \frac {\log n}{n}
\geq \frac {\log 2}{2} + \int_3^{N+1} \frac {\log t}{t} \, dt
\geq \frac {1}{2} \log^2 N - k_3
\end {equation}
with $k_3 = \frac {1}{2} (\log^2 3 - \log 2) = 0.256\ldots$
We have
\begin {eqnarray*}
2 \sum_{A=1}^N \frac {\tau (A)}{A}
& = & 2 \sum_{A=1}^N \sum_{1 \leq m, n : mn = A} \frac {1}{mn}
= 2 \sum_{1 \leq m, n : mn \leq N} \frac {1}{mn}
= 2 \sum_{m=1}^N \frac {1}{m} \sum_{n=1}^{\lfloor N/m \rfloor} \frac {1}{n} \\
& \leq & 2 \sum_{m=1}^N \frac {1}{m} \left( \log \frac {N}{m} + \gamma
   + \frac {1}{2 \left\lfloor \frac {N}{m} \right\rfloor} \right) \\
& \leq & \log^2 N + 4 \gamma \log N + 2 \gamma^2 + \frac {\log N + \gamma}{N}
   + 2 k_3 + \sum_{m=1}^N \frac {1}{m \left\lfloor \frac {N}{m} \right\rfloor} \\
&& \text {by (\ref {eq:harmonic}) and (\ref {eq:sumlog})}
\end {eqnarray*}
Using (\ref {eq:harmonic}), the last term of this sum can be bounded by
\begin {eqnarray*}
\sum_{m=1}^N \frac {1}{m \left\lfloor \frac {N}{m} \right\rfloor}
& \leq & \sum_{m=1}^{\left\lfloor \frac {N}{2} \right\rfloor}
   \frac {1}{m \left( \frac {N}{m} - 1\right)}
   + \sum_{m=\left\lfloor \frac {N}{2} \right\rfloor + 1}^N \frac {1}{m}
\leq 2 \sum_{m=\left\lceil \frac {N}{2} \right\rceil}^N \frac {1}{m}
\leq 2 \log 2 + \frac {1}{N}
\end {eqnarray*}
Combining the inequalities, the height of the class polynomial is bounded from
above by
\begin {equation}
\label {eq:heightwithh}
k_5 h + \pi \sqrt {|D|} \left( \log^2 N + 4 \gamma \log N + \frac {\log N +
\gamma + 1}{N} + k_4 \right)
\end {equation}
with $N = \sqrt {\frac {|D|}{3}}$,
$k_4 = 2 k_3 + 2 \log 2 + 2 \gamma^2 = 2.566\ldots$ and
$k_5 = \log (2 k_2) = 3.011\ldots$

Similarly to the argumentation above one shows that
\begin {equation}
\label {eq:h}
h \leq 2 \sum_{A=1}^N \tau (A) \leq 2 N (\log N + \gamma) + 1,
\end {equation}
which also implies $h \in O \left( \sqrt {|D|} \log |D| \right)$.

Combining (\ref {eq:heightwithh}) and (\ref {eq:h}) yields the bound of
Theorem~\ref {th:height}.

\section {Class group enumeration}
\label {sec:class group}

Computing the class group of $\OD$ by enumerating all reduced primitive quadratic forms of discriminant $D$ is a step of the algorithm that can be neglected in practice. In the new algorithms of Sections~\ref {ssec:multi} and \ref {ssec:agm}, however, its theoretical complexity risks to come close to that of the crucial parts. One has to determine all coprime $[A_i, B_i, C_i]$ such that $|B_i| \leq A_i \leq C_i$, $B_i > 0$ if one of the inequalities is not strict, and $D = B_i^2 - 4 A_i C_i$. These conditions imply that $A_i \leq \sqrt {\frac {|D|}{3}}$.

\subsection {The naïve algorithm}

The following trivial algorithm is implemented in the code described in more detail in
Section~\ref {sec:implementation}: Loop over $A$ and $B$ such that
$0 \leq B \leq A \leq \sqrt {\frac {|D|}{3}}$ and $B \equiv D \pmod 2$ and compute the corresponding $C$, if
it exists. While it can be seen in the figures of Table~\ref {tab:class polynomial}
that this approach takes virtually no time, it carries out $O (|D|)$ arithmetic
operations with integers of bit size in $O (\log |D|)$ and thus has a
complexity of
\[
O (|D| \, M (\log |D|)),
\]
which is smaller than the bound of Corollary~\ref {cor:time} only by logarithmic
factors.

\subsection {Saving one loop}

An asymptotically faster (and again up to logarithmic factors optimal) algorithm
is inspired by the theoretical considerations of Section~\ref {sec:height}.

Looping over $A$ in the interval $1 \leq A \leq \sqrt {\frac {|D|}{3}}$, one
solves the congruence $B = D \bmod 2$ and $\frac {B^2 - D}{4} = 0 \bmod A$, and
keeps the form if $C = \frac {B^2 - D}{4 A} \geq A$. As shown in Section~\ref
{sec:height}, the total number of $B$ to consider is in
$O (\sqrt {|D|} \log |D|)$.

To find the square root of $D$ modulo $A$, the prime factorisation of $A$ is
required, and the most convenient approach is to enumerate the $A$ in factored
form. To this purpose, one needs the primes up to $\sqrt {\frac {|D|}{3}}$,
which can be obtained in time $O \left( \frac {\sqrt {|D|} \log |D|}{\log\log |D|}
\right)$ by \cite {ab04}. All possible $A$ are then computed with an amortised
cost of one multiplication per value in total time $O (\sqrt {|D|} \, M (\log
|D|))$.

For any given $A$ and prime $p | A$, a square root of $D$ modulo $p$ is computed
in $O (\log p \, M (\log p))$ by Cipolla's algorithm \cite {cip03}, lifting to a
root modulo $p^e$ for $p^e || A$ requires an additional $O (e \, M (\log p))$.
For one $A$, this takes
\[
O \left( \sum_{p^e || A} \log (p^e) \, M (\log p) \right)
\subseteq O (\log A \, M (\log A))
\subseteq O (\log |D| \, M (\log |D|)),
\]
or $O \left( \sqrt {|D|} \log |D| \, M (\log |D|) \right)$ for all $A$.

The roots modulo prime powers have to be recombined by the Chinese remainder
theorem to form a candidate $B$. Organising the computations in a tree with the
leaves indexed by the $p^e || A$, one may use the analogues for numbers of the
fast algorithms for polynomials of Section~\ref {ssec:poly}. The tree has
$\omega (A) \in O (\log A)$ leaves and thus a height of $O (\log\log A)$, and the
root contains a number of $O (\log A)$ digits, so that one computation of
Chinese remainders takes
\[
O (\log\log A \, M (\log A)).
\]
This is to be multiplied by the total number of $B$, which yields an overall
complexity of
\[
O (\sqrt {|D|} \log |D| \log\log |D| \, M (\log |D|)).
\]
So Chinese remaindering is the dominant step of this algorithm for class group enumeration, and the first bound of Theorem~\ref {th:classgroup} is proved.

\subsection {Generating the class group by primes}

Another possible approach is to start with a set of prime forms generating the class group. According to \cite[p.~376]{bac90}, under GRH a system of generators is given by the forms having a prime $A$-value bounded by $6 \log^2 |D|$, and these can be computed in time polynomial in $\log |D|$. Let $\pf_1, \pf_2, \ldots$ denote these generators. One may enumerate the powers of $\pf_1$ in the class group by composing and reducing quadratic forms until reaching the order $e_1$ of $\pf_1$ such that $\pf_1^{e_1} = 1$. Next, one computes the powers of $\pf_2$ until $\pf_2^{e_2}$ lies in the subgroup generated by $\pf_1$, constructed in the previous step; then all combinations $\pf_1^{a_1} \pf_2^{a_2}$ with $0 \leq a_i < e_i$ are added to yield the subgroup $\langle \pf_1, \pf_2 \rangle$. One continues with the powers of $\pf_3$ until $\pf_3^{e_3}$ falls into this subgroup, and so forth. If all computed elements are stored in a hash table or simply in a
matrix indexed by $A$ and $B$, looking them up takes negligible time.
The running time of the algorithm is dominated by $O (h)$ computations in the
class group, each of which takes $O (\log \log |D| \, M (\log |D|))$ by
\cite{sch91fas}. This proves the second bound of Theorem~\ref {th:classgroup}.

\subsection {From class groups to $N$-systems}

When working with a class invariant other than $j$ that is not invariant
under $\Sl$, but only under $\Gamma^0 (N)$ for some $N > 0$, the
representatives of the class group need to be normalised to obtain a coherent
set of algebraic conjugates. Such a normalisation is, for instance, given by
an $N$-system as defined in \cite {sch02}. It has the property that
\[
\gcd (A_i, N) = 1 \text { and } B_i \equiv B_1 \pmod {2 N}
\: \text { for } 1 \leq i \leq h
\]
and may be obtained by applying suitable unimodular transformations to the
original $[A_i, B_i, C_i]$. The coefficients of these transformations are
defined modulo $N$, so that for fixed $N$, trying all possibilities
requires a constant number of arithmetic operations per form and increases
the size of the $A_i$, $B_i$ and $C_i$ by a constant factor. Thus, transforming a system of reduced quadratic forms into an $N$-system requires an additional $O (h \, M (\log |D|))$, which is covered by the previous enumeration of the forms. (In practice, one would use a
more intelligent approach to lower the complexity with respect to $N$, cf.
the constructive proof of Proposition~3 in \cite {sch02}.)

\section {Complexity of class polynomial computation}
\label {sec:eval}

We are now able to provide the generic complexity for
class polynomial computation via floating point approximations. Different
approaches to evaluating the class invariants lead to algorithms with
different overall complexities; these are examined below.

Let $n$ be the precision in bits used for the computations,
let $h$ be the class number of $\OD$ and denote by $E (h, n)$ the time
needed to evaluate the class invariant with precision $n$ in $h$ values.
Then the algorithm takes
\begin {itemize}
\item
$O (h^{1 + \varepsilon})$ for enumerating the class group according to
Section~\ref {sec:class group};
\item
$E (h, n)$ for evaluating the class invariant and
\item
$O ((h \log^2 h + \log n) M (n))$ for reconstructing the class polynomial from
its roots according to Algorithm~\ref {alg:poly}.
\end {itemize}

The class group computation is indeed negligible since $E (h, n)$ is at least
of order~$h n$, which is the time needed to write down the conjugates.
So it remains to examine in more detail the quantity $E (h, n)$.

We propose four different algorithms for evaluating modular functions. The
first two of them are well-known, the third one is a novel application of the
techniques of symbolic computation presented in Section~\ref {ssec:poly},
and it already allows to obtain the complexity stated in
Corollary~\ref {cor:time}. The fourth one gains an additional logarithmic
factor for the evaluation phase and yields the slightly more precise
statement of Theorem~\ref {th:time} without
changing the conclusion of Corollary~\ref {cor:time}.

The class polynomial $H_D [f] = \prod_{i=1}^h (X - f (\tau_i))$ is obtained
by evaluating the function $f$ in the $h$ different arguments
$\tau_i = \frac {- B_i + \sqrt D}{2 A_i}$ for an $N$-system $[A_i, B_i,
C_i]$ (see Section~\ref {ssec:poly}), where $N$ depends on $f$ and
is assumed to be fixed. Since $f (z)$ is modular for $\Gamma^0 (N)$, it is
invariant under the translation $z \mapsto z + N$ and admits a Fourier
transform, that is, a Laurent series expansion in the variable
$q^{1/N} = e^{2 \pi i z/N}$.
Thus, the $f (\tau_i)$ may be obtained by first computing the corresponding
$q_i^{1/N}$ and then evaluating the $q$-expansion in these arguments.
By Section~\ref {ssec:multiprec}, the $q_i$ can be computed in time
$O (h \log n \, M (n))$, which will be dominated by
the actual function evaluations.

\subsection {The naïve approach}
\label {ssec:naive}
The straightforward technique for evaluating the modular function
$f = \sum_{\nu = \nu_0}^{\nu_1} c_\nu \left( q^{1/N} \right)^\nu$ consists of
a Horner scheme for the polynomial part
$\sum_{\nu = 0}^{\nu_1 - \nu_0} c_{\nu + \nu_0} \left( q^{1/N} \right)^\nu$
and a multiplication by $\left( q^{1/N} \right)^{\nu_0}$. (Notice that usually
$\nu_0 < 0$.) Its complexity is
$O ((\nu_1 - \nu_0 + \log \nu_0) M (n))$, the $c_\nu$ having been precomputed to
precision $n$. Actually, $\nu_0$ is a constant depending only on $f$;
$\nu_1$, however, depends not only on the $c_\nu$, but also on the desired
precision $n$ and on $\left| q^{1/N} \right|$ and thus on the function argument.

Consider first the classical case of $f$ being $j$, which is invariant under
$\Sl$. We may then assume that the arguments are transformed by a matrix in $\Sl$ into the standard fundamental domain $\Fc$ of (\ref {eq:fundamental domain}) prior to evaluating $j$, so that
$|q| \leq e^{- \pi \sqrt 3}$ is bounded from above by a constant less than $1$.
On the other hand, it is shown in \cite {bp05} that
$0 \leq c_\nu \leq \frac {e^{4 \pi \sqrt \nu}}{\sqrt 2 \, \nu^{3/4}}$
for $\nu \geq 1$.
Thus, it is possible to fix $\nu_1 \in O (n)$ to obtain a precision
of $O (n)$ digits.

The total complexity of evaluating $j$ at $h$ values then becomes
\[
O (h n \, M (n)),
\]
or
\[
O \left( |D| \log^3 |D| \, M \left( \sqrt {|D|} \log^2 |D| \right) \right)
\subseteq O \left( |D|^{3/2} \log^{6 + \varepsilon} |D| \right)
\]
with $n \in O \left( \sqrt {|D|} \log^2 |D| \right)$ and
$h \in O \left( \sqrt {|D|} \log |D| \right)$
according to Section~\ref {sec:height}.

Concerning alternative class invariants, unfortunately the
fundamental domain for $\Gamma^0 (N)$ with $N > 1$ contains at least one
rational number (called a {\em cusp}). In such a cusp, $|q| = 1$, and the
$q$-expansion usually diverges;
in a neighbourhood of the cusp, it may converge arbitrarily slowly.
In this case, it is possible to use a different expansion in the
neighbourhood by transporting the cusp to infinity via a matrix in $\Sl$.
We will not pursue this discussion, since the approach of Section~\ref
{ssec:sparse} provides a faster and simpler solution for all currently
used class invariants.

\subsection {Using the sparsity of $\eta$}
\label {ssec:sparse}

Virtually all class invariants suggested in the literature are in some way
derived from Dedekind's $\eta$-function. This is the case for the Weber
functions
$\f (z) = e^{- \pi i / 24} \, \frac {\eta \left( \frac {z+1}{2} \right)}
{\eta (z)}$,
$\f_1 (z) = \frac {\eta \left( \frac {z}{2} \right)}{\eta (z)}$
and $\f_2 (z) = \sqrt 2 \, \frac {\eta (2z)}{\eta (z)}$ already examined in
\cite {web08},
the generalised Weber functions
$\wN (z) = \frac {\eta \left( \frac {z}{N} \right)}{\eta (z)}$
suggested in \cite {em07},
the double $\eta$ quotients
$\wpq (z) =
\frac {\eta \left( \frac {z}{p_1} \right)\eta \left( \frac {z}{p_2} \right)}
{\eta (z) \eta \left( \frac {z}{p_1 p_2} \right)}$
proposed in \cite {es04},
and even for $j$.
In fact, $j$ is most conveniently computed as
$j = \left( \frac {\f_1^{24} + 16}{\f_1^8} \right)^3$.

The definition of $\eta$ in \cite {ded76} is closely related
to the partition generating function:
\[
\eta = q^{1/24} \, \prod_{\nu \geq 1} (1 - q^\nu);
\]
evaluating the product to precision $n$ requires $O (n)$ arithmetic
operations.
An expression better suited for computation is given by Euler's
pentagonal number theorem \cite {eul83}:
\[
\eta = q^{1/24} \, \left( 1 + \sum_{\nu=1}^\infty
(-1)^\nu \left( q^{\nu (3 \nu - 1) / 2} + q^{\nu (3 \nu + 1) / 2} \right)
\right).
\]
Since the occurring exponents are values of quadratic polynomials, the
series is very sparse: To reach an exponent of order $O (n)$,
only $O (\sqrt n)$ terms need to be computed, and this
process can be implemented with $O (\sqrt n)$ multiplications.
Recall that any polynomial of fixed degree can be evaluated in an
arithmetic progression with a constant number of arithmetic operations per
additional value, once the first few values are known; the employed
algorithm relies on iterated differences. In the
special case of $\eta$, the following recursion yields two additional terms
of the series at the expense of four multiplications by recursively
computing $q^\nu$, $q^{2 \nu -1}$, $q^{\nu (3 \nu - 1) / 2}$ and
$q^{\nu (3 \nu + 1) / 2}$ as follows:
\begin {eqnarray*}
q^{\nu} & = & q^{\nu - 1} \cdot q \\
q^{2 \nu - 1} & = & q^{2 (\nu - 1) - 1} \cdot q^2 \\
q^{\nu (3 \nu - 1) / 2} & = &
q^{(\nu - 1) (3 (\nu - 1) + 1) / 2} \cdot q^{2 \nu - 1} \\
q^{\nu (3 \nu + 1) / 2} & = &
q^{\nu (3 \nu - 1) / 2} \cdot q^\nu
\end {eqnarray*}

Besides the sparse and regular series expression, the $\eta$ function has a
second crucial property that makes it well suited for computation: It is a
modular form of weight $1/2$. As such, unlike $j$, it is not invariant
under transformations in $\Sl$. However, its transformation behaviour is
explicitly known (cf. \cite [\S4]{deu58}) and easily computable.
Thus, to obtain $\eta (z)$
for an arbitrary value of $z$, one should first transform $z$ into the
fundamental domain $\Fc$, so that the series can be truncated at an
exponent of order $O (n)$. Then the evaluation of $O (\sqrt n)$
terms of the $\eta$ series for $O (h)$ distinct values (the constant being
at most~$4$ for the $\wpq$ mentioned in the beginning of this section) with
a floating point precision of $O (n)$ digits can be carried out in time
\[
O (h \sqrt n \, M (n)),
\]
or
\[
O \left( |D|^{3/4} \log^2 |D| \, M \left( \sqrt {|D|} \log^2 |D| \right) \right)
\subseteq O \left( |D|^{5/4} \log^{5 + \varepsilon} |D| \right)
\]
as $h \in O \left( \sqrt {|D|} \log |D| \right)$ and
$n \in O \left( \sqrt {|D|} \log^2 |D| \right)$
according to Section~\ref {sec:height}.

It remains, however, to verify that transforming the arguments
into the fundamental domain is dominated by the cost of the series
evaluation. The arguments being roots of an $N$-system $[A_i, B_i, C_i]$
with $A_i, |B_i| \in O (N^2 \sqrt {|D|}) = O (\sqrt {|D|})$,
they may be transformed into $\Fc$ by reducing the quadratic
forms in time $O (h \, M (\log |D|))$ (see \cite [Prop.~5.4.3]{coh93}),
which is negligible. The same holds for arguments
such as $\frac {z+1}{2}$ or $\frac {z}{N}$, corresponding to quadratic
forms whose discriminants have absolute values in $O (|D|)$.

\subsection {Multipoint evaluation}
\label {ssec:multi}

The algorithms of Sections~\ref {ssec:naive} and \ref {ssec:sparse} compute
the values of modular functions one at a time; but for the sake of class
polynomial computation, we need the values in many points. For polynomials,
Algorithm~\ref {alg:multi} provides a fast way of doing exactly this.
And indeed,
from a numerical point of view a class invariant can be seen as a
polynomial via its truncated $q$-expansion. Either one considers the
function directly as done for $j$ in Section~\ref {ssec:naive}, or one
proceeds via $\eta$ as in Section~\ref {ssec:sparse}. In both cases, a
polynomial of degree $O (n)$ has to be evaluated in $O (h)$ points, which
by Algorithm~\ref {alg:multi} can be done in time
\begin {eqnarray*}
&& O \left( (n \log n + h \log^2 h) M (n) \right)
\subseteq O \left( \sqrt {|D|} \log^3 |D| \, M \left( \sqrt {|D|} \log^2 |D| \right) \right) \\
&& \subseteq O \left( |D| \log^{6 + \varepsilon} |D| \right) .
\end {eqnarray*}

\subsection {Newton iterations on the arithmetic-geometric mean}
\label {ssec:agm}

A new approach for evaluating modular functions in single arguments is described
in \cite {dup07}. It is based on the arithmetic-geometric mean and Newton
iterations on a function involving it. The basic algorithm
underlying \cite [Theorem~4]{dup07} computes the modular
function~$k'$, whose square $\lambda$ satisfies
\[
\frac {256 \left( 1 - \lambda + \lambda^2 \right)^3}
      {\left( \lambda (1 - \lambda) \right)^2} = j.
\]
For an argument with imaginary part bounded by a constant and the
precision~$n$ tending to infinity, it has a complexity of
$O (\log n \, M (n))$.

During class polynomial computations, the precision and the imaginary part of
the arguments are tightly coupled, so that this algorithm is not sufficient to
derive the desired complexity result. The modification of
\cite [Theorem~5]{dup07} obtains the same complexity of $O (\log n \, M (n))$
uniformly in the argument. If the imaginary part of the argument is
of the order of the required precision, then the algorithm of Section~\ref
{ssec:sparse} or even the naïve algorithm of Section~\ref {ssec:naive}
already yield the desired result with a constant number of arithmetic
operations. Otherwise, the argument is repeatedly divided by $2$ until its
imaginary part is smaller than a constant, which can be compensated by
iterations of the arithmetic-geometric mean. Then the previous
Newton algorithm converges sufficiently fast.

For other modular functions $f$, one may have the evaluation of $k'$ followed
by Newton iterations on the modular polynomial relating $k'$ and $f$. For fixed
$f$, this phase does not increase the complexity. The approach does not work,
however, for $\eta$, which is a modular form of weight~$1/2$ instead of a
modular function (of weight~$0$). The algorithm of \cite [Section~7.2]{dup07}
computes first $\theta_{00}^2$, a certain modular form of weight $1$, as the
inverse of
the arithmetic-geometric mean of $1$ and $k'$, and then $\eta$ as the twelfth
root of $\lambda (1 - \lambda) \theta_{00}^2 / 16$ by a suitably initialised
Newton process. Again, one obtains a complexity of $O (\log n \, M (n))$ for an
evaluation at precision $n$, uniformly in the argument.

The complexity for evaluating in $h$ arguments then becomes
\[
O (h \log n \, M (n))
\subseteq O \left( \sqrt {|D|} \log^2 |D| \, M \left( \sqrt {|D|} \log^2 |D| \right) \right)
\subseteq O \left( |D| \log^{5 + \varepsilon} |D| \right)
\]
with the estimates of Section~\ref {sec:height} for $h$ and $n$.
Taking into account the time needed to compute the class polynomial
from its roots by Algorithm~\ref {alg:poly},
this proves Theorem~\ref {th:time}.

\section {Implementation}
\label {sec:implementation}

The algorithms of this article have been implemented using \texttt {gmp} \cite
{gmp} with an assembly patch for 64~bit AMD processors \cite {gau05ass},
\texttt {mpfr} \cite {mpfr} and \texttt {mpc} \cite {mpc} for the
multiprecision arithmetic and \texttt {mpfrcx} \cite {mpfrcx}
for the polynomial
operations. Table~\ref {tab:class polynomial} provides running times for class
numbers between $2500$ and $100000$, obtained on an
AMD Opteron 250 with $2.4$~GHz. All timings are given in seconds and rounded to
two significant digits. (The computations for class number 100000 have been
carried out on a 2.2 GHz machine disposing of more memory, and the running times
have been scaled accordingly.)
For each class number, the discriminant with smallest absolute value has been
chosen. Only the algorithms of Sections~\ref {ssec:sparse} to \ref {ssec:agm}
are taken into account; the naïve approach of Section~\ref {ssec:naive} is
clearly inferior to the one exploiting the sparsity of $\eta$. The
chosen class invariant is the double $\eta$-quotient $\w_{3, 13}$, and
its values are obtained by precomputing a table for the values of $\eta$ at the
$h$ reduced quadratic forms.

\begin {table}[hbt]
\begin {center}
\small
\begin {tabular}{cl|r|r|r|r|r|}
&& $5000$ & $10000$ & $20000$ & $40000$ & $100000$ \\
\hline
& $|D|$ & \footnotesize $6961631$ & \footnotesize $23512271$
 & \footnotesize $98016239$ & \footnotesize $357116231$
 & \footnotesize $2093236031$ \\
(1) & precision $n$ (bits) & $9540$ & $20317$ & $45179$ & $96701$
 & $264727$ \\
(2) & height (in base 2) & $8431$ & $18114$ & $40764$ & $87842$
 & $242410$ \\
(3) & $M (n)$ & $7.3$ & $23$ & $75$ & $230$ & $1080$ \\
\hline
(4) & class group &
   $0.1$ & $0.1$ & $0.4$ & $1.3$ & $6.8$ \\
(5) & conjugates from $\eta$ &
   $3.4$ & $21$ & $140$ & $890$ & $10000$ \\
(6) & poly. from roots &
   $13$ & $93$ & $730$ & $5200$ & $120000$ \\
\hline
& \textbf {sparse series} & & & & & \\
(7) & $\eta$ &
   $12$ & $98$ & $900$ & $7700$ & $140000$\\
(8) & \quad of which $q_i$ &
   $3.0$ & $22$ & $170$ & $1300$ & $20000$\\
(9) & total time &
   $28$ & $210$ & $1800$ & $14000$ & $270000$ \\
\hline
& \textbf {multipoint eval.} & & & & & \\
(10) & $\eta$ &
   $93$ & $640$ & $5700$ & $42000$ & aborted \\
(11) & total time &
   $110$ & $750$ & $6500$ & $48000$ & \\
\hline
     & (10) / (7) & $7.8$ & $6.5$ & $6.3$ & $5.5$ & --- \\
\hline
& \textbf {AGM} & & & & & \\
(12) & $\eta$ &
   $32$ & $200$ & $1400$ & $9900$ & $130000$ \\
(13) & total time &
   $48$ & $320$ & $2300$ & $16000$ & $260000$ \\
\hline
     & (12) / (7) & $2.7$ & $2.0$ & $1.6$ & $1.3$ & $0.93$ \\
\hline
\end {tabular}
\end {center}
\caption {\label {tab:class polynomial}Running times}
\end {table}

The first lines of the table provide some general information.
The precision (1) of the floating point computations is obtained by
increasing the estimate of Section~\ref {sec:height} by $1 \%$ to account for
potential rounding errors. As the chosen class invariant is not $j$, a
correction factor
depending on the class invariant needs to be used, see Section~\ref
{sec:height}. This factor is correct only asymptotically, which explains why the
actual height (2) is a bit smaller than the precision estimate. For $j$, the two
are closer to each other. $M (n)$ is measured by computing the first $100000$
successive powers of $\pi + i \gamma$ with Euler's constant $\gamma$.

The second block of lines provides timings for the steps that are
independent of the algorithm used for evaluating the modular function.
As can be seen, the class group computation (4) is completely negligible.
Line (5) corresponds to the
effort of deriving all values of the class invariant from the tabulated $\eta$
values (reduction of quadratic forms, multiplication by $24$-th roots of unity
and computation of the $\eta$ quotients). The computation of the polynomial from
its roots (6) corresponds precisely to Algorithm~\ref {alg:poly}, and provides a
measure for the complexity of the operations with polynomials. For constructing
the largest polynomial of degree $100000$, the polynomial FFT
has been disabled during the last steps and replaced by Toom--Cook
multiplication, since the FFT consumed too much memory; this explains the jump
from $h=40000$.

The third block contains the timings for evaluating $\eta$ in the reduced
quadratic forms (7) using the sparse series representation as described in
Section~\ref {ssec:sparse}, and the total running time for computing the class
polynomials using this technique (9). Line (8) details the time spent in (7)
(and also in (10)) for computing the $q_i$; it essentially measures the complex
exponential.

The fourth block represents the
corresponding results for the multipoint evaluation approach of
Section~\ref {ssec:multi}, and the last block corresponds to the
as\-ymp\-tot\-i\-cal\-ly
fastest evaluation of Section~\ref {ssec:agm}, for which an implementation by
Dupont has been used. As explained in Section~\ref {ssec:agm}, the algorithm
requires to switch to the sparse series evaluation when the imaginary part of
the argument becomes too large. In the implementation, the AGM code is disabled
for an imaginary part larger than $5$.

Comparing first the evaluation of $\eta$ as a sparse series or by multipoint
evaluation, one notices that the asymptotically faster algorithm is about $5$
to $8$ times slower on the examples and that it appears to catch up with growing
class numbers. However, this happens so slowly that one cannot expect it to beat
the algorithm in $O \left( |D|^{1.25 + \varepsilon} \right)$
in the foreseeable future for any tractable instance.

The approach using Newton iterations on the AGM is faster than multipoint
evaluation, but still hardly beats the asymptotically slower algorithm in $O
\left( |D|^{1.25 + \varepsilon} \right)$: The biggest computed
example of class number $100000$ lies just beyond the cross-over point!

One notices that the growth rates of the running times of all algorithms,
instead of behaving like $|D|$ or $|D|^{1.25}$, come closer to
$|D|^{1.4}$ or $|D|^{1.6}$. A small part of this can be explained by the rather
peculiar choice of discriminants. Taking the first one with a given class
number, the precision $n$ is rather large compared to $h$ and $\sqrt {|D|}$,
since there are many forms with small $A$ so that the sum $\sum_{i=1}^h \frac
{1}{A_i}$ becomes comparatively large.

However, the major reason for the faster than predicted growth of computing time
is that for the floating point arithmetic the range for asymptotically fast
algorithms is not yet reached. The threshold for switching to the FFT in \texttt
{gmp} on the test machine is set to about $500000$ bits; so the
examples lie still in the Karatsuba or Toom-Cook range, which accounts for
a growth of $M (n)$ of $n^{\log 3 / \log 2}$ resp. $n^{\log 5 / \log 3}$ instead
of $n$.

A possible improvement of the multipoint evaluation approach consists of
conveniently grouping the arguments. For instance, in the example of class
number $5000$,  the function $\eta$ has to be evaluated in $2501$ arguments
(corresponding to the two ambiguous forms and $2499$ pairs of opposite
non-am\-big\-u\-ous forms). Assuming the worst case of $|q| \approx e^{- \pi
\sqrt 3}$, that is almost reached for the largest values of $A$, one can
approximate $\eta$ by a polynomial of degree $1190$. As multipoint evaluation
should be most efficient when the number of arguments is about half the degree,
it makes sense to perform four evaluations in $625$ resp. $626$ arguments each.
But sorting them by their absolute values, the smaller ones do not actually
require such a high degree approximation of $\eta$. In the example, an
approximation of degree $260$, $551$, $805$ resp. $1190$ is sufficient for the
four chunks of $q$. Then the time used for multipoint evaluation drops from
$93$~s
to
$64$~s. Experimenting with different partitions of the arguments (three resp.
five parts of the same size, parts of different sizes adapted to the degree of
the approximations, etc.) yields similar results, far from competing with the
sparse series evaluation.

Another point to take into account are the space requirements of the
algorithms. When each root of the class polynomial is computed separately, only
$O (h n)$ bits need to be stored, which is linear in the output size. Multipoint
evaluation as described in Section~\ref {ssec:poly}, however, requires that the
tree constructed in Step~1 of Algorithm~\ref {alg:multi} be maintained in
memory, so that the occupied space grows by a logarithmic factor to reach $O (h
n \log h)$. This logarithmic factor could be saved by evaluating in chunks of
$O \left( \frac {h}{\log h} \right)$ arguments, as explained in Lemma~2.1 of
\cite {gs92}.

Anyway, the example class polynomial of degree $100000$ uses over $5$~GB as an
uncompressed text file and is computed in about $3$ days. This shows that the
limiting factor is the memory requirement rather than the running time, as can
be expected from algorithms that have a close to linear complexity with respect
to their output size.

As a final remark, one notices that the algorithms behave numerically
well, even though rounding errors do occur during floating point computations.
For the algorithm of Section~\ref {ssec:sparse}, this can be explained by
the sparsity of the $\eta$ series and the fact that all coefficients are $+1$
and $-1$.
Indeed, if the last few digits
of a term are erroneous, these errors propagate to subsequent terms. However,
the absolute magnitude of such errors decreases rapidly, so that the wrong
digits in later terms actually do not intervene in the additions. (Otherwise
said, the computations may as well be carried out with fixed point numbers, and
indeed a simulation of fixed point arithmetic using floating point numbers of
decreasing precision yields accurate results.)

\section {Comparison to other approaches}

\subsection {Chinese remaindering}

In \cite {alv04}, the authors suggest an approach for directly computing class
polynomials for $j$ modulo a prime $p$. The basic idea is to derive the
polynomial modulo many small primes by enumerating all elliptic curves modulo
these small primes and only retaining those having complex multiplication by
$\OD$. Then a Chinese remainder technique allows to obtain the class polynomial
modulo~$p$.

Unfortunately, it is not sufficient to gather only $O (\log p)$ bits of
information modulo small primes per coefficient of the class polynomial,
although this is the information contained in the final output. In fact, so
many small primes are needed that the class polynomial could be reconstructed
over $\Z$ instead of only modulo $p$. The complexity derived in Section~3.2 of
\cite {alv04} is
\[
O \left( |D|^{3/2} \log^{10} |D| + |D| \log^2 |D| \log p + \sqrt {|D|} \log^2 p
\right),
\]
and already the term depending only on $|D|$ is worse than what is obtained with
the algorithms of Section~\ref {sec:eval}.

Thus, it is asymptotically faster to compute the class polynomial over $\Z$
using floating point approximations and to reduce it modulo $p$ afterwards.

\subsection {$p$-adic algorithms}

Couveignes and Henocq suggest in \cite {ch02} a $p$-adic approach for computing
class polynomials for $j$. The basic idea is to look for a small prime $p$ and
(by exhaustive search) an elliptic curve modulo $p$ with complex multiplication
by $\OD$. This curve is then lifted to the $p$-adic numbers with a high enough
precision so that the class polynomial may be reconstructed. Alternatively, they
show how to work with supersingular curves. The complexity of their approach
is $O (|D|^{1 + \varepsilon})$, where the exact power of the logarithmic factor
has not been worked out. By nature, the algorithm is not affected by rounding
errors, and using the explicit bound of Theorem~\ref {th:height} its output is
certified to be correct.

So the complex and the $p$-adic approach are both essentially linear in the size
of the class polynomials. All variants have been implemented (for more details
on the complex implementation, see \cite {em02}; for the ordinary $p$-adic
algorithm, see \cite {bs04}; for the supersingular one, \cite {lr04}), and all
seem to work reasonably well in practice. The floating point algorithms are easy
to implement with arbitrary class invariants using the results of \cite
{sch02}. The $p$-adic approach with ordinary curves has been made to work with
certain class invariants other than $j$, see \cite {bs04} and
\cite [Chapter~6]{bro06}.
The considerable overhead involved makes it unclear whether it is competitive
with the complex approach.

\subsection* {Acknowledgements}
I am grateful to R\'egis Dupont for letting me use his code to test the fast
evaluation of $\eta$ via the arithmetic-geometric mean.
I thank Karim Belabas for pointing out to me the basic idea of Section~\ref
{sec:height}, Pierrick Gaudry and Fran\c{c}ois Morain for numerous
discussions and Dan Bernstein, Reinier Br\"oker, Kristin Lauter and Marco
Streng for helpful comments.

\bibliographystyle {amsplain}

\providecommand{\bysame}{\leavevmode\hbox to3em{\hrulefill}\thinspace}
\providecommand{\MR}{\relax\ifhmode\unskip\space\fi MR }
\providecommand{\MRhref}[2]{%
  \href{http://www.ams.org/mathscinet-getitem?mr=#1}{#2}
}
\providecommand{\href}[2]{#2}

\end {document}